\documentclass[11pt]{article}
\usepackage{amsmath,amssymb,amsfonts,amsthm,graphicx}
\usepackage{hyperref}
\newtheorem{theorem}{Theorem}[section]
\newtheorem{proposition}[theorem]{Proposition}

\newtheorem{lemma}[theorem]{Lemma}

\newtheorem{definition}[theorem]{Definition}

\numberwithin{equation}{section}
\title{On convergence rate in the Gauss-Kuzmin problem for $\theta$-expansions}
\author{
    Gabriela Ileana Sebe\footnote{e-mail: igsebe@yahoo.com.} \\
    \emph{\small Politehnica University of Bucharest, Faculty of Applied Sciences},\\
    \emph{\small Splaiul Independentei 313, 060042, Bucharest, Romania} and \\
    \emph{\small Institute of Mathematical Statistics and Applied Mathematics}, \\
     \emph{\small Calea 13 Sept. 13, 050711 Bucharest, Romania} \\
    and\\
    Dan Lascu\footnote{e-mail: lascudan@gmail.com.}\nonumber \\
    \emph{\small Mircea cel Batran Naval Academy, 1 Fulgerului, 900218 Constanta,
    Romania} \\
    }
\begin{document}
\maketitle
\thispagestyle{empty}
\begin{abstract}
After providing an overview of $\theta$-expansions introduced by Chakraborty and Rao, we focus on the Gauss-Kuzmin problem for this new transformation. Actually, we complete our study on these expansions by proving a two-dimensional Gauss-Kuzmin theorem.
More exactly, we obtain such a theorem related to the natural extension of the associated measure-dynamical system.
Finally, we derive explicit lower and upper bounds of the error term which provide interesting numerical calculations for the convergence rate involved.
\end{abstract}

\section{Introduction}

Chakraborty and Rao \cite{CR-2003} have introduced the continued fraction expansion of a number in terms of an irrational $\theta \in (0, 1)$. This new expansion of positive reals is called \textit{$\theta$-expansion}. We mention that the case $\theta = 1$ refers to regular continued fraction (RCF) expansions.
The study initiated by Chakraborty and Rao on the analogous transformation of the Gauss map, was completed by Chakraborty and Dasgupta. Actually, in \cite{CD-2004} was identified the absolutely continuous invariant probability measure of this new transformation only in the particular case $\theta^2 = 1/m$, $m \in \mathbb{N}_+$.

It is only recently that Sebe and Lascu \cite{SL-2014} proved the first Gauss-Kuzmin theorem for $\theta$-expansions, applying the method of random systems with complete connections (RSCC) by Iosifescu and Grigorescu \cite{IG-2009}.
Following the treatment in the case of the RCF, the Gauss-Kuzmin problem for this new transformation can be approached in terms of the associated Perron-Frobenius operator under the invariant measure induced by the limit distribution function.
Moreover, using a Wirsing type approach Sebe \cite{Sebe-2017} obtained a near-optimal solution for the Gauss-Kuzmin problem. The strategy was to restrict the domain of the Perron-Frobenius operator to the Banach space of all functions which have a continuous derivative on $[0, \theta]$.

The aim of this paper is to show a two-dimensional Gauss-Kuzmin theorem for $\theta$-expansions. Note that in the literature there are known similar results for other types of expansions (see \cite{DK-1994, I-1997, DK-1999, Sebe-2000, Sebe-2001}).

The paper is organized as follows.
In the next section we gather prerequisites needed to prove our results in sections 3 and 4. More exactly, in Section 3 we obtain a Gauss-Kuzmin theorem related to the natural extension \cite{Nakada-1981} of the measure-dynamical system corresponding to these expansions.
In Section 4 we try to get close to the optimal convergence rate. Here, the characteristic properties of the Perron-Frobenius operator on the Banach space of functions of bounded variations allows us to derive explicit lower and upper bounds of the error term which provide a more refined estimate of the convergence rate involved.
In the last section we conclude by giving numerical calculations.

\section{Prerequisites}

For a fixed $\theta \in (0,1)$, Chakraborty and Rao \cite{CR-2003} showed that any $x \in \left(0, \theta \right)$ can be written in the form
\begin{equation}
x = \frac{1}{\displaystyle a_1\theta
+\frac{1}{\displaystyle a_2\theta
+ \frac{1}{\displaystyle a_3\theta + \ddots} }} := [a_1 \theta, a_2 \theta, a_3 \theta, \ldots], \label{1.1}
\end{equation}
which is called the $\theta$-\textit{expansion} of $x$.
Here $a_n \in \mathbb{N}_+ : = \left\{1, 2, 3, \ldots\right\}$.
Such $a_n$'s are called $\theta$-\textit{expansion digits} and they are obtained using the transformation
\begin{equation}
T_{\theta}: [0,\theta] \to [0,\theta];\quad
T_{\theta}(x):=
\left\{
\begin{array}{ll}
{\displaystyle \frac{1}{x} - \theta \left \lfloor \frac{1}{x \theta} \right\rfloor} &
{\displaystyle \hbox{if } x \in (0, \theta],}\\
\\
0 & \hbox{if } x=0.
\end{array}
\right. \label{1.2}
\end{equation}
Thus, if we define the {\it quantized index map} $\eta : [0, \theta] \to {\mathbb N}:=\{0,1,2,\ldots\}$ by
\begin{equation}
\eta(x) := \left\{\begin{array}{ll}
\left\lfloor \displaystyle \frac{1}{x \theta}\right\rfloor  & \hbox{if }  x \neq 0, \\
\\
\infty & \hbox{if }  x = 0
\end{array} \right. \label{1.3}
\end{equation}
then the sequence $(a_{n})_{n \in \mathbb N_+}$ in (\ref{1.1}) is obtained as follows:
\begin{equation}
a_n(x) = \eta(T_{\theta}^{n-1}(x)), \quad n \geq 1, \label{1.4}
\end{equation}
with $T_{\theta}^0 (x) = x$.

\begin{figure}[h!]
  \centering
  \includegraphics[width=1\textwidth]{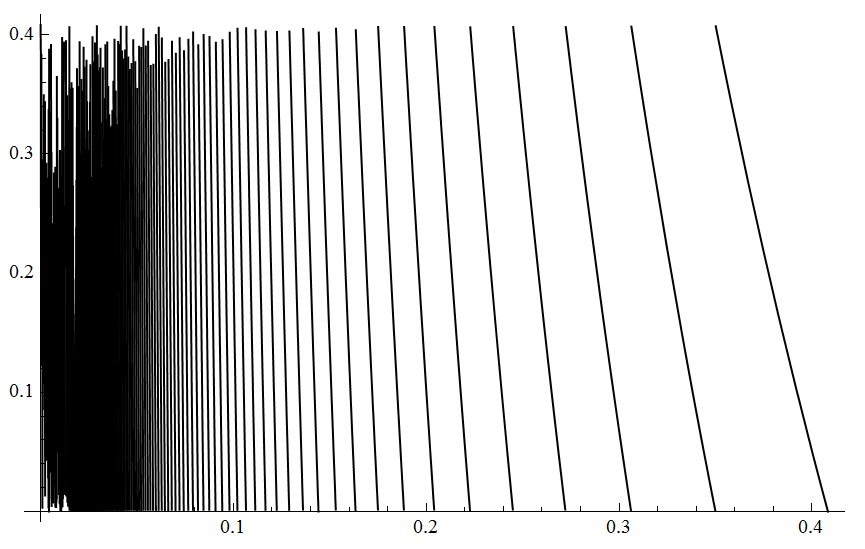}
  \caption{The transformation $T_{\theta}$ for $\theta=\frac{1}{\sqrt{6}}$}\label{fig1}
\end{figure}

This new expansion of positive reals, different from the regular continued fraction expansion, was also studied in
\cite{CD-2004, LN-2017, Sebe-2017, SL-2014}.

In \cite{CR-2003} it was shown that $\theta$-expansions are convergent.
To this end, define real functions $p_n(x)$ and $q_n(x)$, for $n \in \mathbb{N}_+$, by
\begin{eqnarray}
p_n(x) &:=& a_n(x) \theta p_{n-1}(x) + p_{n-2}(x), \quad \label{1.5} \\
q_n(x) &:=& a_n(x) \theta q_{n-1}(x) + q_{n-2}(x), \quad \label{1.6}
\end{eqnarray}
with $p_{-1}(x) := 1$, $p_0(x) := 0$, $q_{-1}(x) := 0$ and $q_{0}(x) := 1$.
It follows that $p_n(x) / q_n(x) = [a_1 \theta, a_2 \theta, \ldots, a_n \theta]$ which is called the \textit{$n$-th order convergent} of $x \in [0, \theta]$. One easily shows that for any $x \in [0, \theta]$ it follows
\begin{equation}
\left|x - \frac{p_n(x)}{q_n(x)}\right| \leq \frac{1}{q_n(x)q_{n+1}(x)} \leq \frac{1}{\left(1+\theta^2\right)^{2\lfloor n/2 \rfloor}\theta^2}, \quad n \in \mathbb{N}_+. \label{1.7}
\end{equation}
Then $p_n(x) / q_n(x) \rightarrow x$, $n \rightarrow \infty$.

In \cite{CR-2003}, Chakraborty and Rao showed that for $\theta^2 = 1/m$, $m \in \mathbb{N}_+$,
$T_{\theta}$ is ergodic with respect to the measure $\gamma_{\theta}$ defined by
\begin{equation}\label{1.8}
\gamma_{\theta} (A) := \frac{1}{\log \left(1+\theta^{2}\right)}
\int_{A} \frac{\theta dx}{1 + \theta x}, \quad A \in {\mathcal{B}}_{[0,\theta]}.
\end{equation}
Let us note that $\gamma_{\theta}$ is $T_{\theta}$-invariant, that is, $\gamma_{\theta}\left(T_{\theta}^{-1}(A)\right)=\gamma_{\theta}(A)$ for any $A \in {\mathcal{B}}_{[0, \theta]}$.
Therefore, $(a_n)_{n \in \mathbb{N}_+}$ is a strictly stationary sequence on $(I, {\mathcal{B}_{[0, \theta]}}, \gamma_{\theta})$ and $a_n \geq m$ for any $m \in \mathbb{N}_+$.

Put $\mathbb{N}_m := \{m, m+1, \ldots\}$, $m \in \mathbb{N}_+$.
For any $n \in \mathbb{N}_+$ and $i^{(n)}=(i_1, \ldots, i_n) \in \mathbb{N}_m^n$ we will say that
\begin{equation}
I \left(i^{(n)}\right) = \{x \in [0, \theta]:  a_k(x) = i_k \mbox{ for } k =1, \ldots, n \} \label{1.9}
\end{equation}
is the $n$\textit{-th order cylinder} and make the convention that $I \left(i^{(0)}\right) = [0, \theta]$.
For example, for any $i \in \mathbb{N}_m$ we have
\begin{equation}
I\left(i\right) = \left\{x \in [0, \theta]: a_1(x) = i \right\} = \left( \frac{1}{(i+1)\theta}, \frac{1}{i \theta} \right). \label{1.10}
\end{equation}

\subsection{Natural extension, extended random variables and Perron-Frobenius operators}

Let $m \in \mathbb{N}_+$ and an irrational $\theta \in (0,1)$ with $\theta^2 = 1/m$.
In this section, we introduce the natural extension $\overline{T}_{\theta}$ of $T_{\theta}$ in (\ref{1.2})
and its extended random variables according to Chap. 1.3 of \cite{IK-2002}, and  we consider the Perron-Frobenius operator of $T_{\theta}$.

\subsubsection{Natural extension}

Let $\left([0, \theta],\mathcal{B}_{[0, \theta]},T_{\theta}\right)$ be as in Section 2.1.
Define $(u_{i})_{i \in \mathbb{N}_m}$ by
\begin{equation}
u_{i}:[0, \theta] \rightarrow [0, \theta]; \quad
u_{i}(x) := \frac{1}{x + i \theta}, \quad x \in [0, \theta]. \label{1.11}
\end{equation}
For each $i \in \mathbb{N}_m$, $u_{i}$ is a right inverse of $T_{\theta}$, that is,
\begin{equation} \label{1.12}
(T_{\theta} \circ u_{i})(x) = x, \quad \mbox{for any } x \in [0, {\theta}].
\end{equation}
Furthermore, if $\eta(x)=i$, then $(u_{i} \circ T_{\theta})(x)=x$ where $\eta$ is as in (\ref{1.3}).
\begin{definition} \label{NatExt}
The \textit{natural extension} $\left([0, {\theta}]^2,{\mathcal B}^2_{[0, {\theta}]},\overline{T}_{\theta}\right)$
of $\left([0, {\theta}],{\mathcal B}_{[0, {\theta}]},T_{\theta}\right)$
is the transformation $\overline{T}_{\theta}$ of the square space
$\left([0, {\theta}]^2,{\mathcal B}^2_{[0, {\theta}]}\right):=\left([0, {\theta}],{\mathcal B}_{[0, {\theta}]}\right) \times \left([0, {\theta}],{\mathcal B}_{[0, {\theta}]}\right)$ defined as follows \cite{Nakada-1981}:
\begin{equation}
\overline{T}_{\theta}:[0, {\theta}]^2 \rightarrow [0, {\theta}]^2;\quad
\overline{T}_{\theta}(x,y) := \left( T_{\theta}(x), \,u_{\eta(x)}(y) \right), \, (x, y) \in [0, {\theta}]^2. \label{1.13}
\end{equation}
\end{definition}

From (\ref{1.12}), we see that $\overline{T}_{\theta}$ is bijective on $[0, {\theta}]^2$ with the inverse
\begin{equation}
(\overline{T}_{\theta})^{-1}(x, y)
= \left(u_{\eta(y)}(x), \, T_{\theta}(y)\right), \quad (x, y) \in [0, {\theta}]^2. \label{1.14}
\end{equation}
Iterations of (\ref{1.13}) and (\ref{1.14}) are given as follows for each
$n \geq 2$:
\begin{eqnarray}
\left( \overline{T}_{\theta} \right)^n(x, y) &=&
(\,T^n_{\theta}(x), \,[x_n{\theta}, x_{n-1}{\theta}, \ldots, x_2{\theta},\, x_1{\theta}+y ] \,), \label{1.15}
\\
\left(\overline{T}_{\theta}\right)^{-n}(x, y) &=&
(\,[y_n{\theta}, y_{n-1}{\theta}, \ldots, y_2{\theta}, \,y_1{\theta}+x ],\, T^{n}_{\theta}(y) \,) \label{1.16}
\end{eqnarray}
where
$x_i:=\eta(T_{\theta}^{i-1}(x))$ and $y_i:=\eta(T_{\theta}^{i-1}(y))$ for $i=1,\ldots,n$.

For $\gamma_{\theta}$ in (\ref{1.8}), define its \textit{extended measure} $\overline{\gamma}_{\theta}$
on $\left([0, {\theta}]^2,{\mathcal{B}}^2_{[0, {\theta}]}\right)$ as
\begin{equation}
\overline{\gamma}_{\theta} (B) := \frac{1}{\log(1 + \theta^2)} \int\!\!\!\int_{B} \frac{\mathrm{d}x\mathrm{d}y}{(1 + xy)^2},
\quad B \in {\mathcal{B}}^2_{[0, \theta]}. \label{1.17}
\end{equation}
Then
$
\overline{\gamma}_{\theta}(A \times [0, \theta])
= \overline{\gamma}_{\theta}([0, \theta] \times A) = \gamma_{\theta}(A)
$
for any $A \in {\mathcal{B}}_{[0, \theta]}$.

The measure $\overline{\gamma}_{\theta}$ is preserved by $\overline{T}_{\theta}$ \cite{SL-2014}, i.e.,
$\overline{\gamma}_{\theta} ((\overline{T}_{\theta})^{-1}(B))= \overline{\gamma}_{\theta} (B)$ for any $B \in {\mathcal{B}}^2_{{\theta}}$.
Since $\overline{T}_{\theta}$ is invertible on $[0, {\theta}]^2$, the last equation is equivalent to
\begin{equation}
\overline{\gamma}_{\theta} \left(\overline{T}_{\theta}(B)\right) = \overline{\gamma}_{\theta} (B),
\quad \mbox{for any } B \in {\mathcal{B}}^2_{[0, {\theta}]}. \label{1.18}
\end{equation}

\subsubsection{Extended random variables}

Define the projection $E:[0, {\theta}]^2 \rightarrow [0, {\theta}]$ by $E(x,y):=x$.
With respect to $\overline{T}_{\theta}$ in (\ref{1.2}), define \textit{extended incomplete quotients} $\overline{a}_l(x,y)$,
$l \in \mathbb{Z}:=\{\ldots, -2,-1,0,1,2, \ldots\}$ at $(x, y) \in [0, {\theta}]^2$ by
\begin{equation}
\overline{a}_{l}(x, y) := (\eta \circ E)(\,(\overline{T}_{\theta})^{l-1} (x, y) \,), \quad l \in \mathbb{Z}. \label{1.19}
\end{equation}
Remark that $\overline{a}_{l}(x, y)$ in (\ref{1.19}) is also well-defined for $l \leq 0$ because $\overline{T}_{\theta}$ is invertible.
For any $n \in \mathbb{N}_+$ and $(x, y) \in [0, {\theta}]^2$, by (\ref{1.15}) and (\ref{1.16}), we have
\begin{equation}
\overline{a}_n(x, y) = x_n, \quad
\overline{a}_0(x, y) = y_1, \quad
\overline{a}_{-n}(x, y) = y_{n+1},
\end{equation} \label{1.20}
where we use notations in  (\ref{1.15}) and (\ref{1.16}).

Since $\overline{\gamma}_{\theta}$ is preserved by $\overline{T}_{\theta}$, the doubly infinite sequence
$(\overline{a}_l(x,y))_{l \in \mathbb{Z}}$ is strictly stationary (i.e., its distribution is invariant under a shift of the indices) under $\overline{\gamma}_{\theta}$.
The stochastic property of $(\overline{a}_l)_{l \in \mathbb{Z}}$ follows from the fact that
\begin{equation}
\overline{\gamma}_{\theta} ( [0, x] \times [0, \theta] \,|
\,\overline{a}_0, \overline{a}_{-1}, \ldots )
=\frac{(a \theta +1)x}{(ax + 1)\theta} \quad \overline{\gamma}_{\theta} \mbox{-}\mathrm{a.s.}, \label{1.21}
\end{equation}
for any $x \in [0, \theta]$, where $a:= [\overline{a}_0 \theta, \overline{a}_{-1} \theta, \ldots]$
with $\overline{a}_{l}:=\overline{a}_l(x,y)$ for $l \in {\mathbb Z}$ and $(x,y) \in [0, \theta]^{2}$.
Hence
\begin{equation}
\overline{\gamma}_{\theta} (\left.\overline{a}_1 = i\right| \overline{a}_0, \overline{a}_{-1}, \ldots) = P_{i}(a) \quad \overline{\gamma}_{\theta} \mbox{-}\mathrm{a.s.}, \label{1.22}
\end{equation}
where
\begin{equation}
P_{i}(x) :=\frac{x \theta + 1}{(x + i\theta)(x + (i+1)\theta )}. \label{1.23}
\end{equation}
The strict stationarity of $\left(\overline{a}_l\right)_{l \in \mathbb{Z}}$, under $\overline{\gamma}_{\theta}$ implies that
\begin{equation}
\overline{\gamma}_{\theta}
(\left.\overline{a}_{l+1} = i\, \right|\, \overline{a}_l,
\overline{a}_{l-1}, \ldots)
= P_{i}(a) \quad \overline{\gamma}_{\theta}
\mbox{-}\mathrm{a.s.} \label{1.24}
\end{equation}
for any $i \in \mathbb{N}_m$ and $l \in \mathbb{Z}$.

The last equation emphasizes that $\left(\overline{a}_l\right)_{l \in \mathbb{Z}}$ is an infinite-order-chain in the theory of dependence with complete connections (see {\cite{IG-2009}}, Section 5.5).

Motivated by (\ref{1.21}), we shall consider the one-parameter family $\{\gamma_{\theta, a}:a \in [0, \theta]\}$ of (conditional) probability measures  on $([0, \theta],{\mathcal{B}}_{[0, \theta]})$ defined by their distribution functions
\begin{equation}
\gamma_{\theta,a} ([0, x]) := \frac{(a \theta +1)x}{(ax + 1)\theta}, \quad x \in [0, \theta], \ a \in [0, \theta]. \label{1.25}
\end{equation}
Note that $\gamma_{\theta,0} = \lambda_{\theta}$.

Let $a_{n}$'s be as in (\ref{1.4}). For each $a \in [0, \theta]$, define $(s_{n,a})_{n \in {\mathbb N}_+}$ by
\begin{equation*}
s_{0,a} := a,\quad s_{n,a} := \frac{1}{a_n \theta + s_{n-1,a}}, \quad n \in \mathbb{N}_+.
\end{equation*}
Then we have
\begin{equation*}
s_{1,a} = \frac{1}{a_1 \theta+a}, \quad  s_{n,a} = \left[a_n \theta, \ldots, a_2 \theta, a_1 \theta+a\right], \, n \geq 2.
\end{equation*}
Note that
\begin{equation*}
\gamma_{\theta,a}\left(A| a_1, \ldots, a_n\right) = \gamma_{\theta,s_{n,a}}\left(T^{n}_{\theta}(A)\right),
\end{equation*}
for all $a \in [0, \theta]$, $A \in \sigma(a_{n+1}, \ldots)$ and $n \in \mathbb{N}_+$.
In particular, it follows that for any $a \in [0, \theta]$
\begin{equation}
\gamma_{\theta,a} \left( \left.T^n_{\theta} < x\right| a_1, \ldots, a_n \right) = \frac{(s_{n,a} \theta +1)x}{(s_{n,a}x + 1)\theta} \label{1.26}
\end{equation}
for any $x \in [0, \theta]$, $n \in \mathbb{N}_+$.
\subsubsection{Perron-Frobenius operator of $T_{\theta}$}

Let $([0, {\theta}],{\mathcal B}_{{\theta}},\gamma_{\theta}, T_{\theta})$ be as in Section 2.1. and let
$L^1([0, \theta],\gamma_{\theta}):=\{f: [0, \theta] \rightarrow \mathbb{C} : \int^{\theta}_{0} |f |\mathrm{d}\gamma_{\theta} < \infty \}$.
The \textit{Perron-Frobenius operator} of $([0, {\theta}],{\mathcal B}_{[0, {\theta}]},\gamma_{\theta}, T_{\theta})$ is defined as the bounded linear operator $U$ on the Banach space $L^1([0, {\theta}],\gamma_{\theta})$ such that the following holds \cite{SL-2014}:
\begin{equation}
Uf(x) = \sum_{i \in \mathbb{N}_m}P_{i}(x)\,f(u_{i}(x)), \quad f \in L^1([0, {\theta}],\gamma_{\theta}) \label{1.27}
\end{equation}
where $P_{i}$ and $u_{i}$ are as in (\ref{1.23}) and (\ref{1.11}), respectively.

For a function $f: [0, \theta] \rightarrow {\mathbb C}$, define the \textit{variation} $\mathrm{var}_{A}f$ of $f$ on a subset $A$ of $[0, \theta]$ by
\begin{equation}
{\rm var}_A f := \sup \sum^{k-1}_{i=1} |f(t_{i+1}) - f(t_{i})|, \label{1.28}
\end{equation}
where the supremum being taken over $t_1 < \cdots < t_k$, $t_i \in A$, $i =1,\ldots,k$ and $k \geq 2$.
We write simply $\mathrm{var} f$ for $\mathrm{var}_{[0, \theta]} f$.
Let $BV([0, \theta]):=\{f:[0, \theta] \rightarrow {\mathbb C}: {\rm var}\,f<\infty\}$ and let
$L^{\infty}([0, \theta])$ denote the collection of all bounded measurable functions $f:[0, \theta] \rightarrow \mathbb{C}$.
It is known that $BV([0, \theta]) \subset L^{\infty}([0, \theta]) \subset L^1([0, \theta], \gamma_{\theta})$.
In \cite{LN-2017} we showed that
\begin{equation} \label{1.29}
\mathrm{var}\,(Uf) \leq \frac{1}{m+1} \mathrm{var } f
\end{equation}
where $U$ is as in (\ref{1.27}) and $f \in BV([0, \theta])$ is a real-valued function.

If $f \in L^{\infty}([0, \theta])$ define the linear functional $U_{\infty}$ by
\begin{equation}
U_{\infty}:L^{\infty}([0, \theta])\to {\mathbb C};\quad U_{\infty} f = \int^{\theta}_{0}f(x)\,\gamma_{\theta}(\mathrm{d} x). \label{2.30'}
\end{equation}
Then we have
\begin{equation} \label{1.30}
U_{\infty} U^n f = U_{\infty} f \quad \mbox{for any } n \in \mathbb{N}_+.
\end{equation}
\begin{proposition} \label{prop.1.2}
For any $f \in BV([0, \theta])$ and for all $n \in \mathbb{N}_+$ we have
\begin{eqnarray}
 \mathrm{var}\,U^n f &\leq& \frac{1}{(m+1)^n} \cdot \mathrm{var } f, \label{1.31} \\
 \left| U^n f - U_{\infty} f \right| &\leq& \frac{1}{(m+1)^n} \cdot \mathrm{var } f. \label{1.32}
\end{eqnarray}
\end{proposition}
\noindent \textbf{Proof.}
Note that for any $f \in BV([0, \theta])$ and $u \in [0, \theta]$ we have
\begin{eqnarray*}
|f(u)| -  \left| \int^{\theta}_{0} f(x) \gamma_{\theta} (\mathrm{dx}) \right| &\leq& \left| f(u)- \int^{\theta}_{0} f(x) \gamma_{\theta} (\mathrm{dx}) \right| \\
&=&  \left| \int^{\theta}_{0}(f(u) - f(x)) \gamma_{\theta} (\mathrm{dx})\right| \leq \mathrm{var } f,
\end{eqnarray*}
whence
\begin{equation}\label{1.33}
|f| \leq  \left| \int^{\theta}_{0} f(x) \gamma_{\theta} (\mathrm{dx}) \right| + \mathrm{var } f, \quad f \in BV([0, \theta]).
\end{equation}
Finally, (\ref{1.30}) and (\ref{1.33}) imply that
\begin{equation*}
\left| U^n f - U_{\infty} f \right| \leq  \mathrm{var } \left( U^n f - U_{\infty} f \right) = \mathrm{var } \, U^n f.
\end{equation*}
for all $n \in \mathbb{N}$ and $f \in BV([0, \theta])$, which leads to (\ref{1.32}).

\hfill $\Box$

By induction with respect to $n \in \mathbb{N}$ we get
\begin{equation*}
U^n f(x) = \sum_{i_1, \ldots, i_n \in \mathbb{N}_m} P_{i_1 \ldots i_n} (x) f(u_{i_n \ldots i_1}(x)), \quad x \in [0, \theta]
\end{equation*}
where
\begin{eqnarray}
  u_{i_n \ldots i_1}   &=& u_{i_n} \circ \ldots \circ u_{i_1} \label{1.34'}\\
  P_{i_1 \ldots i_n} (x) &=& P_{i_1} (x) P_{i_2} (u_{i_1} (x)) \ldots P_{i_n} (u_{i_{n-1} \ldots i_1} (x)), \quad n \geq 2. \label{1.35}
\end{eqnarray}
Here the functions $u_{i}$ and $P_i$ are defined in (\ref{1.11}) and (\ref{1.23}), respectively, for all $i \in \mathbb{N}_m$.

Putting
\begin{equation*}
\frac{p_n(i_1, \ldots, i_n)}{q_n(i_1, \ldots, i_n)} = [i_1 \theta, \ldots, i_n \theta], \quad n \in \mathbb{N}_+,
\end{equation*}
for arbitrary indeterminates $i_1, \ldots, i_n$, we get
\begin{eqnarray} \label{1.36}
P_{i_1\ldots i_n} (a) &=& \frac{1+a \theta}{q_{n-1}(i_2, \ldots, i_n)(a+i_1 \theta) + p_{n-1}(i_2, \ldots, i_n)} \nonumber \\
                      &\times& \frac{1}{q_{n}(i_2, \ldots, i_n, m)(a+i_1 \theta) + p_{n}(i_2, \ldots, i_n, m)}
\end{eqnarray}
for all $n \geq 2$, $i_n \in \mathbb{N}_m$, and $a \in [0, \theta]$.

\section{Gauss-Kuzmin theorem related to the natural extension}
In this section a Gauss-Kuzmin theorem for $([0, \theta]^2, \mathcal{B}^2_{[0, \theta]}, \overline{\gamma}_{\theta}, \overline{T}_{\theta})$ is given.
First we give a modified version of the Gauss-Kuzmin theorem for $T_{\theta}$ proved in \cite{LN-2017}.
Then we show some important results used in the proof of the main theorem.
\subsection{Gauss-Kuzmin theorems for $T_{\theta}$ and $\overline{T}_{\theta}$}

\begin{theorem} \label{G-K1} (A Gauss-Kuzmin theorem for $T_{\theta}$)
Let $([0, \theta], \gamma_{\theta}, T_{\theta})$ as in (\ref{1.2}) and (\ref{1.8}).
If $\lambda_{\theta}$ is the Lebesgue measure on $[0, \theta]$,
then there exists a constant $0< q < \theta$ such that for any $A \in \mathcal{B}_{[0, \theta]}$ we have
\begin{equation}\label{2.1}
\left| \lambda_{\theta}\left(T_{\theta}^{-n}(A)\right) - \gamma_{\theta}(A) \right| < C \lambda_{\theta}(A) \mathcal {O}(q^n)
\end{equation}
where $C$ is an universal constant.
\end{theorem}
\noindent \textbf{Proof.}
In \cite{SL-2014} (Prop.14(ii)) we show that
$\mu \left((T_{\theta})^{-n}(A)\right) = \int_{A} U^nf(x) \mathrm{d} \gamma_{\theta}(x)$
where $\mu$ is a probability measure on $\left([0, \theta],{\mathcal{B}}_{[0, \theta]}\right)$
absolutely continuous with respect to the Lebesgue measure $\lambda_{\theta}$,
and $f(x):= (\log(1+\theta^2)) \frac{1+\theta x}{\theta^2}h(x)$ with $h := d\mu / d \lambda_{\theta}$ a.e. in $[0, \theta]$.
In the special case $\mu = \lambda_{\theta}$ we obviously have
\begin{equation}\label{3.2}
\lambda_{\theta} \left(T_{\theta}^{-n}(A)\right) = \frac{\theta}{\log\left(1+\theta^2\right)}\int_{A} \frac{U^nf(x)}{1+\theta x} \mathrm{d}x
\end{equation}
with $f(x):= (\log(1+\theta^2)) \frac{1+\theta x}{\theta^2}$, $x \in [0, \theta]$.
Thus, from (\ref{2.30'}) we have that $U_{\infty} f= 1$. Therefore,
\begin{equation}\label{3.3}
\gamma_{\theta}(A) =  \frac{\theta}{\log\left(1+\theta^2\right)}\int_{A} \frac{U_{\infty}f}{1+\theta x} \mathrm{d}x.
\end{equation}
Using \cite{SL-2014} (Prop.27) it follows that there exist two positive constants $q<\theta$ and $K$ such that
\begin{equation}
\left\|U^nf-U_{\infty}f\right\|_L \leq K q^n \left\|f\right\|_L, \quad f \in L([0, \theta]), n \in \mathbb{N}_+ \label{3.4}
\end{equation}
where $L([0, \theta])$ denote the Banach space of all complex-valued Lipschitz continuous functions on $[0, \theta]$ with the following norm:
\begin{equation*}
\left\| f \right\|_L := \sup_{x \in [0, \theta]} \left|f(x)\right| + \sup_{x' \neq x''} \frac{|f(x') - f(x'')|}{|x' - x''|}. \label{3.5}
\end{equation*}
Therefore
\begin{eqnarray*}
\left| \lambda_{\theta}\left(T_{\theta}^{-n}(A)\right) - \gamma_{\theta}(A) \right|
&\leq& \frac{\theta}{\log\left(1+\theta^2\right)} \int_{A} \frac{\left|U^nf(x)-U_{\infty}f\right|}{1+\theta x} \mathrm{d}x \\
&<& K q^n \left\|f\right\|_L \frac{\theta}{\log\left(1+\theta^2\right)} \int_{A} \frac{1}{1+\theta x} \mathrm{d}x \\
&<& C q^n \gamma_{\theta}(A)
\end{eqnarray*}
and since
\begin{equation*}
\gamma_{\theta}(A) \leq \frac{\theta}{\log\left(1+\theta^2\right)} \lambda_{\theta}(A), \quad A \in \mathcal{B}_{[0, \theta]}
\end{equation*}
then the proof is complete.
\hfill $\Box$

In \cite{LN-2017} the proof of Gauss-Kuzmin theorem is based on the Gauss-Kuzmin-type equation which in this case is
\begin{equation}
F_{n+1} (x) = \sum_{i \geq m }\left\{F_{n}\left(\frac{1}{i \theta}\right) - F_{n}\left(\frac{1}{i \theta + x}\right)\right\} \label{2.2}
\end{equation}
where the functions $(F_n)_{n \in \mathbb{N}}$ are defined for $x \in [0,\theta]$ by
\begin{equation}
F_0 (x) := \lambda_{\theta} ([0,x]) = \frac{x}{\theta}, \quad F_n (x) := \lambda_{\theta} (T_{\theta}^n \leq x), \, n \in \mathbb{N}_+. \label{2.3}
\end{equation}

The measure $\gamma_{\theta}$ defined in (\ref{1.8}) is an eigenfunction of (\ref{2.2}), namely,
if we put $F_{n}(x) = \log (1+\theta x)$, $x \in [0, \theta]$, we obtain $F_{n+1}(x) = \log (1+\theta x)$.
The factor $1/(\log(1+\theta^2))$ is a normalizing constant.

We give now the main theorem of this section. First, define $\Delta_{x,y}=[0,x] \times [0,y]$ for any $x,y \in [0, \theta]$,
and the functions $(\overline{F}_{n})_{n \in \mathbb{N}_+}$ on $[0, \theta]^2$ by
\begin{equation} \label{2.4}
\overline{F}_{n}(x,y) := \overline{\lambda}_{\theta} \left( \left(\overline{T}_{\theta}\right)^n(x,y) \in \Delta_{x,y} \right),
\end{equation}
where $\overline{\lambda}_{\theta}$ is the Lebesgue measure on $[0, \theta]^2$.

\begin{theorem} \label{G-K2} (A Gauss-Kuzmin theorem for $\overline{T}_{\theta}$)
For every $n \geq 2$ and $(x,y) \in [0, \theta]^2$ one has
\begin{equation}\label{2.5}
\overline{F}_{n}(x,y) = \frac{\log(1+xy)}{\log(1+\theta^2)} + {\mathcal O} \left(q^n\right)
\end{equation}
with $0<q<\theta$.
\end{theorem}

\subsection{Necessary results}

In this subsection, we give necessary results used to prove the Gauss-Kuzmin theorem for $\overline{T}_{\theta}$.
Like in the one-dimensional case, we need the Gauss-Kuzmin-type equation associated with the functions $(\overline{F}_{n})_{n \in \mathbb{N}_+}$ defined in (\ref{2.4}).
Thus, for any $0 < y \leq \theta$, put $\ell_1 := \eta(y)$, where $\eta$ is as in (\ref{1.3}).
Then $\left(\overline{T}_{\theta}\right)^{n+1}(x,y) \in \Delta_{x,y}$ is equivalent to
\begin{equation*}
\left(\overline{T}_{\theta}\right)^{n} \in \left\{ \bigcup_{i \geq \ell_1+1} \left[ \frac{1}{i \theta +x}, \frac{1}{i \theta} \right] \times [0,\theta] \right\} \cup
\left\{ \left[\frac{1}{\ell_1 \theta +x}, \frac{1}{\ell_1 \theta} \right]  \times \left[ \frac{1}{y} - \ell_1 \theta, \theta \right] \right\}.
\end{equation*}
From this and (\ref{2.4}) we get the Gauss-Kuzmin-type equation on $[0, \theta]^2$:
\begin{eqnarray} \label{2.6}
\overline{F}_{n+1}(x,y) &=& \sum_{i \geq \ell_1} \left\{ \overline{F}_{n}\left(\frac{1}{i \theta}, \theta\right) - \overline{F}_{n}\left(\frac{1}{i \theta +x}, \theta\right) \right\} \nonumber \\
&-& \left\{
\overline{F}_{n}\left(\frac{1}{\ell_1 \theta},\frac{1}{y} - \ell_1 \theta\right) -
\overline{F}_{n}\left(\frac{1}{\ell_1 \theta +x},\frac{1}{y}-\ell_1 \theta\right)
\right\}.
\end{eqnarray}
A straightforward calculation shows that the measure $\overline{\gamma}_{\theta}$ defined in (\ref{1.17}) is an eigenfunction of (\ref{2.6}),
namely, if we put $\overline{F}_{n}(x,y) = \log(1+xy)$, $x,y \in [0, \theta]$,
we obtain $\overline{F}_{n+1}(x,y) = \log(1+xy)$.
\begin{lemma} \label{lema2.3}
Let $n \in \mathbb{N}$, $n \geq 2$ and let $y \in [0, \theta] \cap \mathbb{Q}$ with $y = [\ell_1 \theta, \ldots, \ell_d \theta]$, $\ell_1, \ldots, \ell_d \in \mathbb{N}_m$, $\ell_d \geq m+1$, where $d \leq \lfloor n/2 \rfloor$. Then for every $x, x^* \in [0,\theta)$ with $x^* < x$,
\begin{equation*}
\left| \overline{F}_n(x,y) - \overline{F}_n(x^*,y) - \frac{1}{\log(1+\theta^2)} \log \left( \frac{1+xy}{1+x^*y} \right) \right| <
C \overline{\lambda}_{\theta}(\Delta_{x,y} \setminus \Delta_{x^*,y}) q^{n-d}
\end{equation*}
where $C$ is an universal constant and $q$ is from Theorem \ref{G-K1}.
\end{lemma}
\noindent \textbf{Proof.}
Let $y_0=y$, $y_i:=[\ell_{i+1} \theta, \ldots, \ell_d \theta]$, $i=1, \ldots, d$, with $y_d=0$. Then $y_1=\frac{1}{y}- \ell_1 \theta$.
Applying (\ref{2.6}) one gets
\begin{eqnarray}
\overline{F}_{n}(x,y) - \overline{F}_{n}(x^*,y)
&=& \sum_{i \geq \ell_1} \left\{ \overline{F}_{n-1}\left(\frac{1}{i \theta + x^*}, \theta\right) - \overline{F}_{n-1}\left(\frac{1}{i \theta +x}, \theta\right) \right\} \nonumber \\
&+& \left( \overline{F}_{n-1}\left(\frac{1}{\ell_1 \theta +x}, y_1\right) - \overline{F}_{n-1}\left(\frac{1}{\ell_1 \theta + x^*}, y_1 \right) \right). \label{2.7}
\end{eqnarray}
Now for each $B \in {\mathcal B}^2_{[0, \theta]}$ one has
\begin{equation} \label{2.8}
\frac{1}{(1+\theta^2)} \frac{1}{\log \left(1+\theta^2\right)} \overline{\lambda}_{\theta}(B) \leq \overline{\gamma}_{\theta}(B)
\leq \frac{1}{\log \left(1+\theta^2\right)} \overline{\lambda}_{\theta}(B).
\end{equation}
Now from (\ref{1.10}) and (\ref{2.8}), it follows that:
\begin{eqnarray} \label{2.9}
&&\sum_{i \geq \ell_1} {\lambda}_{\theta} \left( \left[ \frac{1}{i \theta + x}, \frac{1}{i \theta +x^*} \right] \right) =
\sum_{i \geq \ell_1} \frac{1}{\theta}\left( \frac{1}{i \theta +x^*} - \frac{1}{i \theta +x} \right) \nonumber \\
&&=\sum_{i \geq \ell_1} \overline{\lambda}_{\theta} \left\{ \left(\frac{1}{i \theta + x}, \frac{1}{i \theta +x^*}\right) \times [0,\theta]\right\} \nonumber \\
&&\leq (1+\theta^2) \log \left(1+\theta^2\right) \sum_{i \geq \ell_1} \overline{\gamma}_{\theta} \left\{\left(\frac{1}{i \theta + x}, \frac{1}{i \theta +x^*}\right) \times [0,\theta]\right\} \nonumber\\
&&= (1+\theta^2) \log \left(1+\theta^2\right) \sum_{i \geq \ell_1} \overline{\gamma}_{\theta} \left\{ (x^*, x) \times I(i) \right\} \nonumber\\
&&\leq (1+\theta^2) \log \left(1+\theta^2\right) \frac{1}{\log \left(1+\theta^2\right)} \sum_{i \geq \ell_1} \overline{\lambda}_{\theta} \left\{ (x^*, x) \times I(i) \right\}  \nonumber \\
&&= (1+\theta^2) \sum_{i \geq \ell_1} \overline{\lambda}_{\theta} \left\{ (x^*, x) \times \left(\frac{1}{(i+1)\theta}, \frac{1}{i \theta}\right) \right\} \nonumber \\
&&= (1+\theta^2) \frac{x-x^*}{\theta} \sum_{i \geq \ell_1} \lambda_{\theta} \left\{ \left(\frac{1}{(i+1)\theta}, \frac{1}{i \theta}\right) \right\} \nonumber \\
&&= (1+\theta^2) \frac{x-x^*}{\theta} \frac{1}{\theta} \sum_{i \geq \ell_1} \left(\frac{1}{i\theta} - \frac{1}{(i+1) \theta}\right) \nonumber \\
&&\leq (1+\theta^2) \frac{x-x^*}{\theta^2} \frac{1}{\ell_1 \theta} \leq (1+\theta^2) \frac{x-x^*}{\theta} (m+1) \frac{y}{\theta} \nonumber \\
&&= \frac{(1+\theta^2)^2}{\theta^2} \overline{\lambda}_{\theta} \left(\Delta_{x,y} \setminus \Delta_{x^*,y}\right).
\end{eqnarray}

For every $2 \leq k \leq d$, a similar analysis leads to
\begin{eqnarray} \label{2.10}
&& \sum_{i \geq \ell_k} \lambda_{\theta} \left\{ [i\theta, \ell_{k-1}\theta, \ldots, \ell_{1}\theta+x], [i\theta, \ell_{k-1} \theta, \ldots, \ell_{1} \theta +x^*] \right\} = \nonumber \\
&& \frac{1}{\theta}\sum_{i \geq \ell_k} \left| [i\theta, \ell_{k-1}\theta, \ldots, \ell_{1}\theta+x^*] - [i\theta, \ell_{k-1}\theta, \ldots, \ell_{1}\theta+x] \right| = \nonumber \\
&& \sum_{i \geq \ell_k} \overline{\lambda}_{\theta} \left\{ [i\theta, \ell_{k-1}\theta, \ldots, \ell_{1}\theta+x], [i\theta, \ell_{k-1} \theta, \ldots, \ell_{1} \theta +x^*] \times [0, \theta] \right\} \leq \nonumber \\
&& (1+\theta^2) \log \left(1+\theta^2\right) \times  \nonumber \\
&&\sum_{i \geq \ell_k} \overline{\gamma}_{\theta} \left\{ [i\theta, \ell_{k-1}\theta, \ldots, \ell_{1}\theta+x], [i\theta, \ell_{k-1} \theta, \ldots, \ell_{1} \theta +x^*] \times [0, \theta] \right\} = \nonumber \\
&& (1+\theta^2) \log \left(1+\theta^2\right) \sum_{i \geq \ell_k} \overline{\gamma}_{\theta} \left\{ \left(x^*, x\right) \times I(\ell_1, \ldots, \ell_{k-1}) \right\} \leq \nonumber \\
&& \frac{(1+\theta^2) \log \left(1+\theta^2\right)}{\log \left(1+\theta^2\right)} \sum_{i \geq \ell_k} \overline{\lambda}_{\theta} \left\{ \left(x^*, x\right) \times I(\ell_1, \ldots, \ell_{k-1}) \right\} \leq \nonumber \\
&& \frac{(1+\theta^2)^2}{\theta^2} \overline{\lambda}_{\theta} \left(\Delta_{x,y} \setminus \Delta_{x^*,y}\right).
\end{eqnarray}
Since $\overline{F}_n(x,\theta) = F_n(x)$, from Theorem \ref{G-K1} it follows that
\begin{eqnarray} \label{2.11}
&&\sum_{i \geq \ell_1}
\left\{
\overline{F}_{n-1}\left( \frac{1}{i \theta +x^*},\theta\right) - \overline{F}_{n-1}\left( \frac{1}{i \theta +x},\theta\right)
\right\} = \nonumber \\
&&\sum_{i \geq \ell_1}
\left\{
{F}_{n-1}\left( \frac{1}{i \theta +x^*}\right) - {F}_{n-1}\left( \frac{1}{i \theta +x}\right)
\right\} = \nonumber \\
&&\sum_{i \geq \ell_1}
\left\{
\gamma_{\theta} \left( \left[\frac{1}{i\theta +x},\frac{1}{i\theta+x^*} \right]\right) +
\lambda_{\theta}\left( \left[\frac{1}{i\theta +x},\frac{1}{i\theta+x^*} \right]\right) {\mathcal O} \left(q^{n-1}\right)
\right\}= \nonumber \\
&& \frac{1}{\log\left(1+\theta^2\right)} \log{\frac{\ell_1\theta+x}{\ell_1 \theta+x^*}} +
\frac{(1+\theta^2)^2}{\theta^2} \overline{\lambda}_{\theta} \left(\Delta_{x,y} \setminus \Delta_{x^*,y}\right) {\mathcal O} \left(q^{n-1}\right).
\end{eqnarray}
Now from (\ref{2.7}), (\ref{2.11}), we have:
\begin{eqnarray*}
\overline{F}_{n}(x,y) - \overline{F}_{n}(x^*,y)&=&
\frac{1}{\log\left(1+\theta^2\right)} \log{\frac{\ell_1\theta+x}{\ell_1 \theta+x^*}}  \nonumber \\
&+&\frac{(1+\theta^2)^2}{\theta^2} \overline{\lambda}_{\theta} \left(\Delta_{x,y} \setminus \Delta_{x^*,y}\right) {\mathcal O} \left(q^{n-1}\right) \nonumber \\
&+&\left( \overline{F}_{n-1}\left(\frac{1}{\ell_1 \theta +x}, y_1\right) - \overline{F}_{n-1}\left(\frac{1}{\ell_1 \theta + x^*}, y_1 \right) \right).
\end{eqnarray*}
Applying (\ref{2.7}) and Theorem \ref{G-K1} again, we have:
\begin{eqnarray*}
&&\overline{F}_{n-1}\left(\frac{1}{\ell_1 \theta +x},y_1\right) - \overline{F}_{n-1}\left(\frac{1}{\ell_1 \theta +x^*},y_1\right)= \\
&&\sum_{i \geq \ell_2}
\left\{
\overline{F}_{n-2}\left( \frac{1}{i \theta + \frac{1}{\ell_1 \theta +x^*}},\theta\right) - \overline{F}_{n-2}\left( \frac{1}{i \theta +\frac{1}{\ell_1 \theta +x}},\theta\right)
\right\} \\
&&+\overline{F}_{n-2}\left( \frac{1}{\ell_2 \theta + \frac{1}{\ell_1 \theta +x^*}},y_2\right) - \overline{F}_{n-2}\left( \frac{1}{\ell_2 \theta +\frac{1}{\ell_1 \theta +x}},y_2\right) \\
&&= \frac{1}{\log\left(1+\theta^2\right)} \log\frac{\ell_2 \theta+[\ell_1 \theta +x]}{\ell_2 \theta+[\ell_1 \theta +x^*]} \\
&&+\frac{(1+\theta^2)^2}{\theta^2} \overline{\lambda}_{\theta} \left(\Delta_{x,y} \setminus \Delta_{x^*,y}\right) {\mathcal O} \left(q^{n-2}\right)\\
&&+\overline{F}_{n-2}\left( \frac{1}{\ell_2 \theta + \frac{1}{\ell_1 \theta +x^*}},y_2\right) - \overline{F}_{n-2}\left( \frac{1}{\ell_2 \theta +\frac{1}{\ell_1 \theta +x}},y_2\right).
\end{eqnarray*}
Applying (\ref{2.7}) and Theorem \ref{G-K1} $d$-times and tacking into account that $y_d = 0$, we get
\begin{eqnarray*}
&&\overline{F}_{n}(x,y) - \overline{F}_{n}(x^*,y)= \frac{1}{\log\left(1+\theta^2\right)} \times \\
&& \log  \left( \frac{\ell_1 \theta +x}{\ell_1 \theta +x^*} \cdot \frac{[\ell_1 \theta +x]+\ell_2 \theta}{[\ell_1 \theta +x^*]+\ell_2 \theta} \cdot \frac{[\ell_{d-1} \theta, \ldots, \ell_2 \theta, \ell_1 \theta +x]+\ell_d \theta}{[\ell_{d-1} \theta, \ldots, \ell_2 \theta, \ell_1 \theta +x^*]+\ell_d \theta} \right) + \\
&&+\frac{(1+\theta^2)^2}{\theta^2} \overline{\lambda}_{\theta} \left(\Delta_{x,y} \setminus \Delta_{x^*,y}\right) \left( {\mathcal O} \left(q^{n-1}\right) + \ldots+ {\mathcal O} \left(q^{n-d}\right) \right).
\end{eqnarray*}
If $p_d$ and $q_d$ are as in (\ref{1.5}) and (\ref{1.6}) with $a_1=\ell_1 +\frac{x}{\theta}$ and $a_i=\ell_i$, $i=2,\ldots,d$, then
$\ell_1 \theta +x = \displaystyle\frac{q_1}{q_0}$, $\ell_2 \theta + [\ell_1 \theta +x] = \displaystyle\frac{q_2}{q_1}$,
$\ell_d \theta +[\ell_{d-1} \theta, \ldots, \ell_2 \theta, \ell_1 \theta +x] =\displaystyle\frac{q_{d}}{q_{d-1}}$.

Let $p^*_d$ and $q^*_d$ are as in (\ref{1.5}) and (\ref{1.6}), with $a_1=\ell_1 +\frac{x^*}{\theta}$ and $a_i=\ell_i$, $i=2,\ldots,d$.
Note that $p_d = p^*_d$.
Thus we have
\begin{eqnarray*}
&&\frac{\ell_1 \theta +x}{\ell_1 \theta +x^*} \cdot \frac{[\ell_1 \theta +x]+\ell_2 \theta}{[\ell_1 \theta +x^*]+\ell_2 \theta} \cdot \frac{[\ell_{d-1} \theta, \ldots, \ell_2 \theta, \ell_1 \theta +x]+\ell_d \theta}{[\ell_{d-1} \theta, \ldots, \ell_2 \theta, \ell_1 \theta +x^*]+\ell_d \theta} = \frac{q_d}{q^*_d} = \frac{p^*_d}{q^*_d} \frac{q_d}{p_d} \\
&&= \frac{x+\ell_1 \theta+[\ell_2 \theta, \ldots, \ell_d \theta]}{x^*+\ell_1 \theta+[\ell_2 \theta, \ldots, \ell_d \theta]}=
\frac{x+\frac{1}{y}}{x^*+\frac{1}{y}} = \frac{1+xy}{1+x^*y}.
\end{eqnarray*}
Therefore,
\begin{eqnarray*}
\overline{F}_{n}(x,y) - \overline{F}_{n}(x^*,y)=
\frac{1}{\log\left(1+\theta^2\right)} \log \left( \frac{1+xy}{1+x^*y} \right)
+\frac{(1+\theta^2)^2}{\theta^2} \overline{\lambda}_{\theta} \left(\Delta_{x,y} \setminus \Delta_{x^*,y}\right) {\mathcal O} \left(q^{n-d}\right)
\end{eqnarray*}
which completes the proof.
\hfill $\Box$

\subsection{Proof of Theorem \ref{G-K2}}

Let $(x, y) \in [0, \theta]^2$, $n \geq 2$ and $y \notin\mathbb{Q}$.
Since $\Delta_{x,p_d/q_d} \subset \Delta_{x,y}$ and
$\overline{F}_n(x,y)=\overline{\lambda}_{\theta} \left(\left(\overline{T}_{\theta}\right)^{-n}(\Delta_{x,y})\right)$, from (\ref{2.8}), (\ref{1.7}) and the fact that $\overline{T}_{\theta}$ is $\overline{\gamma}_{\theta}$-invariant, we find that
\begin{eqnarray} \label{2.12}
&&\overline{F}_{n}(x,y) - \overline{F}_{n}(x,\frac{p_d}{q_d}) =
\overline{\lambda}_{\theta}\left(\left(\overline{T}_{\theta}\right)^{-n}(\Delta_{x,y}) \setminus \left(\overline{T}_{\theta}\right)^{-n}(\Delta_{x,p_d/q_d}) \right) \nonumber \\
&&\leq \left(1+\theta^2\right) \log \left( 1+\theta^2 \right) \overline{\gamma}_{\theta}\left(\left(\overline{T}_{\theta}\right)^{-n}(\Delta_{x,y}) \setminus \left(\overline{T}_{\theta}\right)^{-n}(\Delta_{x,p_d/q_d}) \right) \nonumber \\
&&= \left(1+\theta^2\right) \log \left( 1+\theta^2 \right) \overline{\gamma}_{\theta} \left(\left(\overline{T}_{\theta}\right)^{-n}\left(\Delta_{x,y} \setminus \Delta_{x,p_d/q_d}\right) \right) \nonumber\\
&&\leq \left(1+\theta^2\right) \log \left( 1+\theta^2 \right) \frac{1}{\log \left( 1+\theta^2 \right)} \overline{\lambda}_{\theta} \left( [0,x] \times \left[ \frac{p_d}{q_d}, y\right] \right) \nonumber \\
&& = \left(1+\theta^2\right) \frac{x}{\theta} \frac{1}{\theta} \left|y - \frac{p_d}{q_d}\right|
\leq \frac{1+\theta^2}{\theta^2} \frac{x}{\left(1+\theta^2\right)^{2\lfloor n/2 \rfloor}\theta^2}.
\end{eqnarray}
Since for every fixed $x \in [0,\theta]$ the function $y \mapsto \log \left(1+xy\right)$ is a differentiable on $[0, \theta]$,
by the \textit{Mean Value Theorem} we have
\begin{eqnarray} \label{2.13}
\left| \log \left(1+xy\right) - \log \left(1+x\frac{p_d}{q_d}\right) \right| = \left| y-\frac{p_d}{q_d} \right| \cdot
\left| \frac{x}{1+x \xi } \right| \leq \frac{x}{\left(1+\theta^2\right)^{2\lfloor n/2 \rfloor}\theta^2}
\end{eqnarray}
where $p_d/q_d \leq \xi \leq y$.

Finally, from Lemma \ref{lema2.3}, (\ref{2.12}) and (\ref{2.13}), and since $\overline{F}_n(0, \frac{p_d}{q_d})=0$, we have:
\begin{eqnarray*}
&&\left| \overline{F}_{n}(x,y) - \frac{\log(1+xy)}{\log\left(1+\theta^2\right)} \right| \leq \left| \overline{F}_{n}(x,y) - \overline{F}_{n}\left(x,\frac{p_d}{q_d}\right)\right|  \\
&&+ \left| \overline{F}_n\left(x,\frac{p_d}{q_d}\right) - \overline{F}_n\left(0,\frac{p_d}{q_d}\right) -
\frac{1}{\log\left(1+\theta^2\right)} \log \left(1+x\frac{p_d}{q_d}\right) \right| \\
&&+ \frac{1}{\log\left(1+\theta^2\right)} \left| \log \left(1+xy\right) - \log \left(1+x\frac{p_d}{q_d}\right) \right| \\
&&\leq \frac{1+\theta^2}{\theta^2} \frac{x}{\left(1+\theta^2\right)^{2\lfloor n/2 \rfloor}\theta^2} + C q^{n-d}+ \frac{x}{\left(1+\theta^2\right)^{2\lfloor n/2 \rfloor}\theta^2}
\end{eqnarray*}
which completes the proof.
\hfill $\Box$

\section{A two-dimensional Gauss-Kuzmin theorem}

In this section we shall estimate the error term
\begin{equation*}
e_{n,a} (x,y) = \gamma_{\theta,a} \left( T^n_{\theta} \in [0,x], s_{n,a} \in [0,y] \right)  - \frac{\log(1+xy)}{\log\left(1+\theta^2\right)}
\end{equation*}
for any $a \in [0, \theta]$, $x, y \in [0, \theta]$ and $n \in \mathbb{N}_+$.

In the main result of this section, Theorem \ref{th.4.4}, we shall derive lower and upper bounds (not depending on $a \in [0, \theta]$) of the supremum
\begin{equation}\label{4.1}
\sup_{x, y \in [0, \theta]} |e_{n,a} (x, y)|, \quad a \in [0, \theta],
\end{equation}
which provide a more refined estimate of the convergence rate involved.
First, we obtain a lower bound for the following approximation error.

\begin{theorem} \label{Th.4.1}
For any $a \in [0, \theta]$ and $n \in \mathbb{N}_+$ we have
\begin{equation*}
\frac{1}{2} P_{m(n)}(\theta) \leq \sup_{y \in [0, \theta]} \left|\gamma_{\theta,a} \left( s_{n,a} \leq y \right) - \gamma_{\theta} \left([0,y] \right) \right|.
\end{equation*}
\end{theorem}
\noindent \textbf{Proof.}
The continuity of the function $y \rightarrow \gamma_{\theta} \left([0,y] \right)$, $y \in [0, \theta]$ and the equation
\begin{equation*}
\lim_{h\searrow 0} \gamma_{\theta, a} \left(s_{n,a} \leq y-h \right) = \gamma_{\theta, a} \left(s_{n,a} \leq y \right)
\end{equation*}
imply that
\begin{equation*}
\sup_{y \in [0, \theta]} \left|\gamma_{\theta,a} \left( s_{n,a} \leq y \right) - \gamma_{\theta} \left([0,y] \right) \right| =
\sup_{y \in [0, \theta]} \left|\gamma_{\theta,a} \left( s_{n,a} < y \right) - \gamma_{\theta} \left([0,y] \right) \right|
\end{equation*}
for any $a \in [0, \theta]$ and $n \in \mathbb{N}_+$. For any $s \in [0, \theta]$ we then have
\begin{eqnarray*}
  \gamma_{\theta,a} (s_{n,a}=s) &=& \gamma_{\theta,a} \left( s_{n,a} \leq s \right) -  \gamma_{\theta} \left([0,s] \right) \nonumber \\
  &-& \left(\gamma_{\theta,a} \left( s_{n,a} < s \right) - \gamma_{\theta} \left([0,s] \right) \right) \nonumber \\
  &\leq& \sup_{y \in [0,\theta]} \left|\gamma_{\theta,a} \left( s_{n,a} \leq y \right) - \gamma_{\theta} \left([0,y] \right) \right| \nonumber \\
  &+& \sup_{y \in [0, \theta]} \left|\gamma_{\theta,a} \left( s_{n,a} < y \right) - \gamma_{\theta} \left([0,y] \right) \right| \nonumber \\
  &=& 2 \sup_{y \in [0, \theta]} \left|\gamma_{\theta,a} \left( s_{n,a} \leq y \right) - \gamma_{\theta} \left( [0,y] \right) \right|.
\end{eqnarray*}
Hence
\begin{equation*}
\sup_{y \in [0, \theta]} \left|\gamma_{\theta,a} \left( s_{n,a} \leq y \right) - \gamma_{\theta} \left( [0,y] \right) \right| \geq
\frac{1}{2} \sup_{s \in [0, \theta]} \gamma_{\theta,a} \left( s_{n,a} = s \right),
\end{equation*}
for any $a \in [0, \theta]$ and $n \in \mathbb{N}_+$.
Next, using (\ref{1.35}) we have
\begin{equation*}
\gamma_{\theta,a} \left( s_{n,a} = [i_n \theta, \ldots, i_2 \theta, i_1 \theta +a] \right) = P_{i_1\ldots i_n}(a), \ n\geq 2,
\end{equation*}
\begin{equation*}
\gamma_{\theta,a} \left( s_{1,a} = \frac{1}{i_1 \theta +a} \right) = P_{i_1}(a)
\end{equation*}
for any $a \in [0, \theta]$ and $i_1,\ldots,i_n \in \mathbb{N}_m$.
By (\ref{1.36}) we have
\begin{equation*}
\sup_{s \in [0, \theta]} \gamma_{\theta,a} \left( s_{n,a} = s \right) = P_{m(n)}(a), \quad a \in [0, \theta]
\end{equation*}
where we write $m(n)$ for $(i_1, \ldots, i_n)$ with $i_1=\ldots=i_n=m$, $n \in \mathbb{N}_+$, $m \in \mathbb{N}_+$.

By the same equation we have
\begin{eqnarray*}
P_{m(n)} (a) &=& \frac{1+a \theta}{q_{n-1}(\underbrace{m, \ldots, m}_{(n-1) \ times})(a+m \theta) + p_{n-1}(\underbrace{m, \ldots, m}_{(n-1) \ times})} \nonumber \\
                      &\times& \frac{1}{q_{n}(\underbrace{m, \ldots, m, m}_{n \ times})(a+m \theta) + p_{n}(\underbrace{m, \ldots, m, m}_{n \ times})}.
\end{eqnarray*}
It is easy to see that $P_{m(n)} (\cdot)$ is a decreasing function. Therefore
\begin{equation*}
\sup_{s \in [0, \theta]} \gamma_{\theta,a} \left( s_{n,a} = s \right) \geq P_{m(n)}(\theta), \quad n \in \mathbb{N}_+
\end{equation*}
for any $a \in [0, \theta]$.
\hfill $\Box$

\begin{theorem} \label{th.4.2}
(The lower bound)
For any $a \in [0, \theta]$ and $n \in \mathbb{N}_+$ we have
\begin{equation*}
\frac{1}{2} P_{m(n)}(\theta) \leq \sup_{x, y \in [0, \theta]} \left|\gamma_{\theta,a} \left( T^n_{\theta} \in [0,x], s_{n,a} \in [0,y] \right) - \frac{\log (1+xy)}{\log\left(1+\theta^2\right)} \right|.
\end{equation*}
\end{theorem}
\noindent \textbf{Proof.} For any $a \in [0, \theta]$ and $n \in \mathbb{N}_+$, by Theorem \ref{Th.4.1} we have
\begin{eqnarray*}
&\sup_{x,y \in [0, \theta]}& \left|\gamma_{\theta,a} \left( T^n_{\theta} \in [0,x], s_{n,a} \in [0,y] \right) - \frac{\log(1+xy)}{\log\left(1+\theta^2\right)} \right| \nonumber \\
&\geq& \sup_{y \in [0, \theta]} \left|\gamma_{\theta,a} \left( T^n_{\theta} \in [0, \theta], s_{n,a} \in [0,y] \right) -
\frac{\log(1+\theta y)}{\log\left(1+\theta^2\right)} \right| \nonumber \\
&=& \sup_{y \in [0, \theta]} \left|\gamma_{\theta,a} \left( s_{n,a} \in [0,y] \right) - \gamma_{\theta} \left( [0,y] \right) \right| \geq \frac{1}{2} P_{m(n)}(\theta).
\end{eqnarray*}
\hfill $\Box$

In what follows we use the characteristic properties of the transition operator associated with the RSCC underlying $\theta$-expansions. By restricting this operator to the Banach space of functions of bounded variation on $[0, \theta]$, we derive an explicit upper bound for the supremum (\ref{4.1}).

\begin{theorem} \label{th.4.3}
(The upper bound)
For any $a \in [0, \theta]$ and $n \in \mathbb{N}$ we have
\begin{equation*}
\sup_{x, y \in [0, \theta]} \left|\gamma_{\theta,a} \left( T^n_{\theta} \in [0,x], s_{n,a} \in [0,y] \right) - \frac{\log(1+xy)}{\log\left(1+\theta^2\right)} \right| \leq \frac{1}{(m+1)^n}.
\end{equation*}
\end{theorem}
\noindent \textbf{Proof.} Let $F_{n,a}(y) = \gamma_{\theta,a} (s_{n,a} \leq y)$ and $G_{n,a} (y) = F_{n,a}(y) - \gamma_{\theta} ([0,y])$,
$a,y \in [0,\theta]$, $n \in \mathbb{N}$.
As we have noted $U$ is the transition operator of the Markov chain $(s_{n,a})_{n \in \mathbb{N}}$ on $\left([0, \theta], \mathcal{B}_{[0, \theta]}, \gamma_{\theta, a}\right)$ for any $a \in [0, \theta]$.
For any $y \in [0, \theta]$ consider the function $f_y$ defined on $[0, \theta]$ as
\begin{equation*}
f_{y}(a):=
\left\{
\begin{array}{ll}
{1} & { \mbox{if } \, 0 \leq a \leq y,}\\
{0} & \mbox{if } \, y < a \leq \theta.
\end{array}
\right.
\end{equation*}
Hence
\begin{equation*}
U^n f_y (a) = E_a\left( \left. f_y(s_{n,a})\right| s_{0,a} = a \right) = \gamma_{\theta,a} (s_{n,a} \leq y)
\end{equation*}
for all $a,y \in [0, \theta]$, $n \in \mathbb{N}$.
As
\begin{equation*}
U_{\infty} f_y = \int^{\theta}_{0} f_y(a) \gamma_{\theta} (\mathrm{da}) = \gamma_{\theta} ([0,y]), \quad y \in [0, \theta].
\end{equation*}
It follows from Proposition \ref{prop.1.2} that
\begin{eqnarray} \label{4.2}
  |G_{n,a}(y)| &=& \left| \gamma_{\theta,a} (s_{n,a} \leq y) - \gamma_{\theta} ([0,y]) \right| \nonumber \\
               &=& \left| U^n f_y (a) - U_{\infty} f_y \right| \leq \frac{1}{(m+1)^n} \mathrm{var } \, f_y = \frac{1}{(m+1)^n}
\end{eqnarray}
for all $a,y \in [0, \theta]$, $n \in \mathbb{N}$.
By (\ref{1.26}), for all $a,x,y \in [0, \theta]$ and $n \in \mathbb{N}$ we have
\begin{eqnarray*}
  \gamma_{\theta,a} \left( T^n_{\theta} \in [0,x], s_{n,a} \in [0,y] \right) = \int^{y}_{0}  \gamma_{\theta,a} \left( \left.T^n_{\theta} \in [0,x] \right| s_{n,a} = z \right) \mathrm{d F_{n,a}(z)}  \nonumber \\
   = \int^{y}_{0} \frac{(z \theta+1)x}{(zx+1)\theta} \mathrm{d F_{n,a}(z)} = \int^{y}_{0} \frac{(z \theta+1)x}{(zx+1)\theta} {\gamma_{\theta}(\mathrm{dz})} + \int^{y}_{0} \frac{(z\theta+1)x}{(zx+1)\theta} \mathrm{d G_{n,a}(z)} \nonumber \\
   = \frac{\log(1+xy)}{\log\left(1+\theta^2\right)} + \frac{(z \theta+1)x}{(zx+1) \theta} \left.G_{n,a}(z)\right|^{y}_{0}
   - \int^{y}_{0} \frac{x(\theta-x)}{(zx+1)^2 \theta} G_{n,a}(z)\mathrm{dz}.
\end{eqnarray*}
Hence, by (\ref{4.2})
\begin{eqnarray*}
\left|\gamma_{\theta,a} \left( T^n_{\theta} \in [0,x], s_{n,a} \in [0,y] \right) - \frac{\log(1+xy)}{\log\left(1+\theta^2\right)} \right|\\
\leq \frac{1}{(m+1)^n} \left( \frac{(y \theta+1)x}{(xy+1) \theta} - \frac{(\theta-x)xy}{(xy+1)\theta} \right) = \frac{x}{(m+1)^n \theta} \leq \frac{1}{(m+1)^n}
\end{eqnarray*}
for all $a,x,y \in [0, \theta]$ and $n \in \mathbb{N}$.
\hfill $\Box$

Combining Theorem \ref{th.4.2} with Theorem \ref{th.4.3} we obtain Theorem \ref{th.4.4}.
\begin{theorem} \label{th.4.4}
For any $a \in [0, \theta]$ and $n \in \mathbb{N}_+$ we have
\begin{eqnarray*}
\frac{1}{2} P_{m(n)}(\theta) &\leq& \sup_{x, y \in [0, \theta]} \left|\gamma_{\theta,a} \left( T^n_{\theta} \in [0,x], s_{n,a} \in [0,y] \right) - \frac{\log(1+xy)}{\log\left(1+\theta^2\right)} \right| \\
&\leq& \frac{1}{(m+1)^n}.
\end{eqnarray*}
\end{theorem}

\section{Final remarks}
To conclude this paper, we note that
\[
P_{m(n)}(\theta) = \frac{m+1}{q_{n+1}q_{n+2}}, \quad n \in \mathbb{N}_+,
\]
where $q_n = m \theta q_{n-1}+q_{n-2}$, $n \in \mathbb{N}_+$, with $q_{-1}=0$ and $q_{0}=1$.
It is easy to see that
\[
q_n = \frac{\theta}{\sqrt{1+4 \theta^2}}\left[ \left(\frac{1+\sqrt{1+4 \theta^2}}{2 \theta}\right)^{n+1} -  \left(\frac{1-\sqrt{1+4 \theta^2}}{2 \theta}\right)^{n+1}\right].
\]
It should be noted that Theorem \ref{th.4.4} in connection with the limits
\begin{eqnarray*}
  \lim \left(\frac{1}{2} P_{m(n)}(\theta)\right)^{1/n} &=& \frac{2 \theta^2}{1+2 \theta^2 +\sqrt{1+4 \theta^2}},  \\
  \lim \left(\frac{1}{(m+1)^n}\right)^{1/n} &=& \frac{1}{m+1},
\end{eqnarray*}
leads to an estimate of the order of magnitude of the supremum (\ref{4.1}.
Actually, Theorem \ref{th.4.4} implies that the convergence rate is $\mathcal{O}(\alpha^n)$, with
\begin{equation*}
\frac{2 \theta^2}{1+2 \theta^2 +\sqrt{1+4 \theta^2}} \leq \alpha \leq \frac{1}{m+1}.
\end{equation*}

For example, we have

\begin{center}
\begin{tabular}{|l|l|}
  \hline
  $m=1$ & $g^2=0.381966 \leq \alpha \leq 0.50000$\\\hline
  $m=2$ & $0.267949192 \leq \alpha \leq 0.33333\ldots$\\\hline
  $m=3$ & $0.208711948 \leq \alpha \leq 0.25000\ldots$\\\hline
  $m=10$ & $0.083920216 \leq \alpha \leq 0.090909\ldots$\\\hline
  $m=100$ & $0.009804864 \leq \alpha \leq 0.00990099$\\\hline
  $m=1000$ & $0.000998004 \leq \alpha \leq 0.000999$\\\hline
  $m=10000$ & $0.00009998 \leq \alpha \leq 0.00009999$\\
  \hline
\end{tabular}
\end{center}

\begin{figure}[h!]
  \centering
  \includegraphics[width=1\textwidth]{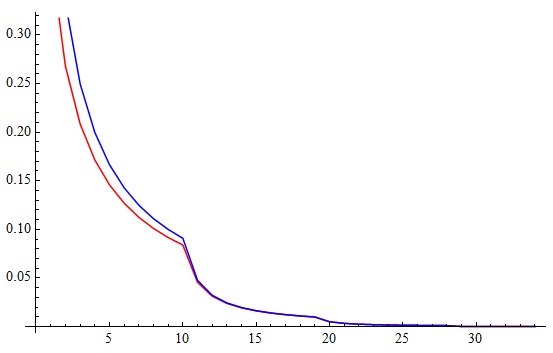}
  \caption{Graphs of lower and upper bounds}\label{fig2}
\end{figure}

\end{document}